\documentclass[11pt]{article}
\usepackage{amsmath, amsthm, amssymb}    
\usepackage{newtxtext}
 
\usepackage[margin=3.14cm]{geometry}
\usepackage{graphicx}			
\usepackage[colorlinks=true, urlcolor=blue, citecolor=black, linkcolor=black]{hyperref}  
\usepackage{gensymb}
\usepackage{float}
\usepackage{graphicx,wrapfig,lipsum}
\usepackage{subcaption}			
\usepackage{xfrac}
\usepackage[symbol]{footmisc}
\usepackage{sectsty}
\sectionfont{\centering}

\usepackage{scrextend}
\usepackage{cite}
\usepackage{url}

\def \r {\tilde{r}}
\def \B {\mathcal{B}}
\def \R {\mathbb{R}}

\renewenvironment{thebibliography}[1]{
  \begin{oldthebibliography}{#1}
    \setlength{\itemsep}{0em}
    \setlength{\parskip}{0em}
}
{
  \end{oldthebibliography}
}

\title{Crocheting Bour's $\B_m$ minimal surfaces}

\author{Hanne Kekkonen\textsuperscript{} 
\vspace{10pt}\\
\textsuperscript{}EEMCS, Delft University of Technology, Netherlands; h.n.kekkonen@tudelft.nl} 

\date{}					

\begin{document}

\maketitle


\begin{abstract}
Minimal surfaces can be though as a mathematical generalisation of surfaces formed by soap films. We consider Bour's minimal surfaces $\B_m$ that are intrinsically surfaces of revolution. We show how to generate crochet patterns for $\B_m$ surfaces using basic trigonometric identities to calculate required arc lengths. Three special cases of $\B_m$ surfaces are considered in more detail, namely Enneper's, Richmond's, and Bour's $\B_3$ surfaces, and we provide exact crochet instructions for the classical Enneper's surface.
\end{abstract}

\section{Introduction}

Mathematics and crochet might not appear the most likely pairing for most people. However, crocheting is an inherently mathematical process. You can create various shapes using stitches with different heights and increasing or decreasing the number of stitches in certain places. Crochet also allows creating many shapes that are very difficult to make with any other technique. One can crochet shapes like the Klein bottle or Seifert surfaces of knots in freestyle, without following exact instructions, since these surfaces do not have a strict shape. A Klein bottle can be short and wide or tall and narrow. But there are many surfaces, for example a sphere or a disc, that have a specific shape. To crochet such models you need crochet instructions and to calculate these instructions you need a good understanding of the underlying mathematical model. 

The idea of knitting or crocheting mathematical or scientific models is not new, though it has not been used very widely. The Scottish chemist Alexander Crum Brown knitted several interlinked surfaces to visualise the ideas presented in his paper "On a Case of Interlacing Surfaces' in the late $19\textsuperscript{th}$ century \cite{Brown1886}. Miles Reid wrote a paper on knitting mathematical surfaces in the 1970s \cite{Reid1971} which inspired several new patterns including a Möbius scarf and the Klein bottle. The crocheted hyperbolic surfaces were introduced by Daina Taimina in 1997 \cite{Henderson2001} and the idea led to a bloom of so called hyperbolic crochet. A few years after the paper on hyperbolic crochet  Hinke Osinga and Bernd Krauskopf described how to crochet an approximation of the Lorenz manifold \cite{Osinga2004}. See also \cite{Taimina2018} for further examples on mathematical crochet.  

Both the hyperbolic plane and the Lorenz manifold require exact crochet instructions. The hyperbolic plane has constant negative Gaussian curvature and so it looks the same at every point. This allows a rather simple pattern that can be worked in rounds and where, after a few set-up rounds, every n$\textsuperscript{th}$ stitch is doubled. 
The Lorenz manifold is a less regular surface and requires a much more complex pattern of stitches. The model is also worked in rounds but unlike the hyperbolic surface it requires detailed instructions on when to add or remove stitches. It takes the full advantage of the versatility of crocheting requiring three different type of stitches which allows different parts of a round to have different heights. In this paper we consider Bour's minimal surfaces $\B_m$ which are 'crochet symmetric' allowing simple crochet instructions (excluding possible intersections) requiring only one type of stitch and with the added or removed stitches evenly spaces across a round.

\section{Minimal surfaces}\label{Sec:MinimalSurfaces}

Minimal surfaces can be considered as a mathematical generalisation of soap film surfaces. If you dip a wire frame into soap water the created surface is optimal in the sense that it minimises the surface area bounded by the frame. 
In mathematics minimal surfaces are defined as surfaces that locally minimise their area. They are allowed to self-intersect and do not have to have a boundary. In general, minimal surfaces do not provide a minimum for the surface area globally. 

The early study of minimal surfaces was rather theoretical. Joseph-Louis Lagrange considered in his 1761 paper \cite{Lagrange1761} the problem of finding a surface $z(x,y)$ of least area given a closed boundary. 
He stated that such minimal surfaces can be locally expressed as the graph of a solution of the minimal surface equation
\begin{align*}
\left(1+z_x^2\right) z_{y y}-2 z_x z_y z_{x y}+\left(1+z_y^2\right) z_{x x}=0
\end{align*}
but the only concrete example he could provide was the plane. Joseph Plateau published his experimental observations on minimal surfaces obtained as soap films bounded by various wires in 1873 \cite{Plateau1873}. His work gave a clear physical interpretation to the problem and helped to spread its study beyond mathematics. As a result the problem of showing the existence of a minimal surface given a boundary became to be known as the Plateau's problem. Only special cases of the problem were solved prior to the early 1930's when the general case was solved independently by Jesse Douglas \cite{Douglas1931} and Tibor Rad\'o \cite{Rado1930}. 

Before the 19$\textsuperscript{th}$ century only two non-trivial minimal surfaces were known; the catenoid discovered by Leonhard Euler in 1744 and the helicoid described by Jean Baptiste Meusnier in 1776. Meusnier also showed that the catenoid satisfies Lagrange's condition and gave a geometric interpretation for the minimal surface equation noting that it is equivalent to the surface having a vanishing mean curvature \cite{Meusnier1785}.
The first general approach for finding minimal surfaces was introduced by Karl Weierstrass in 1866 \cite{Weierstrass1866}. This representation allows construction of minimal surfaces based on an integral of a holomorphic and meromorphic function. The analytical representation of minimal surfaces using holomorphic coordinates was independently formulated by Alfred Enneper a few years prior to the publication of Weierstrass' representation (though published in 1868 \cite{Enneper1868}) which is why this way of constructing minimal surfaces is often called Weierstrass-Enneper representation.

The Weierstrass-Enneper equations of a minimal surface are 
\begin{align*}
x(\zeta) & = \text{Re}\int f(\zeta)(1-g^2(\zeta))d\zeta\\
y(\zeta) & = \text{Re}\int i f(\zeta)(1+g^2(\zeta))d\zeta\\
z(\zeta) & = \text{Re}\int 2f(\zeta)g(\zeta)d\zeta,
\end{align*}
where $f$ and $g$, known as the Weierstrass-Enneper data of the minimal surface, are a holomorphic and a meromorphic function respectively. 
In this paper we consider minimal surfaces that are intrinsically surfaces of revolution determined by the Weierstrass-Enneper data 
\begin{align}\label{eq:EWdata}
f(\zeta) = C\zeta^{m-2}\quad \text{and} \quad g(\zeta)=\zeta,
\end{align}
where $m\in\R$. We denote these surfaces by $\B_m$ after Edmond Bour who studied them in his paper \cite{Bour1862}, see also \cite{Whittemore1917}. 
Choosing $m=0$ and $C=1$ in \eqref{eq:EWdata} produces the catenoid while the choice $m=0$ and $C=i$ leads to the right helicoid.

\begin{figure}[h]
\centering
\includegraphics[width=1\textwidth]{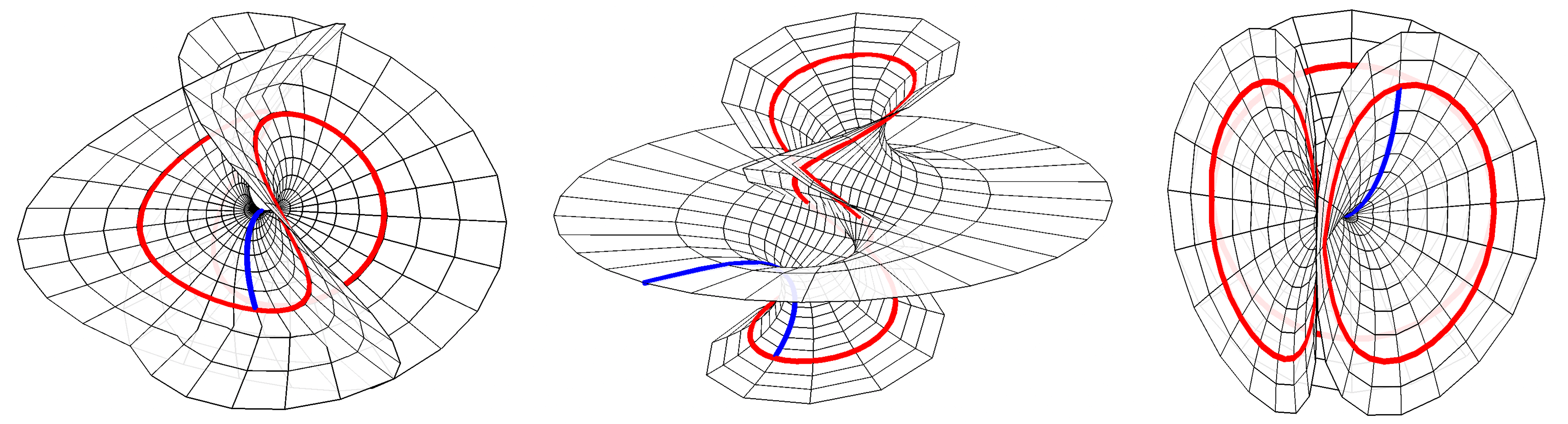}
\vspace{0.5cm}
\caption{To crochet an accurate $\B_m$ surface we need to calculate the circumference of a circle centred at the origin  (red) with an intrinsic radius $R=\ell\cdot H$ (blue) where $H$ is the height of a stitch and $\ell\in\mathbb{N}_+$ is the round.}
\label{Fig:Surfaces}
\end{figure}

If we assume that $C=1$, $m\in\mathbb{R}\setminus\{-1,0,1\}$, and denote $\zeta = r e^{i\theta}$ we can rewrite the above as
\begin{align}\label{eq:GeneralEW}
\begin{split}
x(r,\theta) & = \frac{r^{m-1}}{m-1}\cos((m-1)\theta) - \frac{r^{m+1}}{m+1}\cos((m+1)\theta)\\
y(r,\theta) & = -\frac{r^{m-1}}{m-1}\sin((m-1)\theta) - \frac{r^{m+1}}{m+1}\sin((m+1)\theta)\\
z(r,\theta) & = \frac{2r^m}{m}\cos(m\theta).
\end{split}
\end{align}
Bour's paper was published in 1862, a couple of years before the Weierstrass-Enneper representation, and he defined the surfaces in the coordinates $r$ and $\theta$. Notice that the catenoid and helicoid cannot be presented in the above form since they corresponded to the Weierstrass-Enneper data with $m=0$. 
Choosing $m=2$, $m=1/2$ or $m=3$ we get the classical Enneper's, Richmond's, and Bour's surfaces respectively, see Figure \ref{Fig:Surfaces}. All three are special cases of the Bour's surfaces $\B_m$ but the Enneper surface, and its Weierstrass-Enneper data, were first considered in more detail by Alfred Enneper in 1864, and the Richmond's surface was studied by Herbert William Richmond in 1904. The Bour's $\B_3$ surface is often simply referred to as the Bour's minimal surface. 
For a comprehensive introduction to minimal surfaces see e.g. \cite{Nitsche1989}. 

We will see in the next section that the first fundamental form of the Bour's minimal surfaces is rotationally symmetric meaning that the $\B_m$ are intrinsically surfaces of revolution. This means that the length of a line segment on the surface does not depend on the angle $\theta$ which allows calculating simple crochet instructions where only one type of stitch is required and the  added or removed stitches can be spaced evenly on a round. See Section \ref{Sec:CrochetPractice} and Appendix \ref{Appendix} for crochet instructions for the classical Enneper's surface. 

\section{Calculating crochet instructions for $\B_m$ surfaces}\label{Sec:CrochetTheory}

Before moving to calculating crochet instructions for the Bour's minimal surfaces let us have a look on how to crochet a simple flat disc. Denote the height of a crochet stitch by $H$ and the width by $W$. Then the radius of the crocheted disc will always be $H\cdot\ell$, with $\ell\in\mathbb{N}_+$ denoting the number of the round. Since a circle with radius $R$ drawn on an Euclidean plane has circumference $2\pi R$ the circumference of a crocheted disc on a round $\ell$ should be $2\pi( H\cdot\ell)$ and the number of stitches required to achieve this is $2\pi ( H\cdot\ell)/W$, rounded to the closest integer. If you use a single crochet stitch the width and height of the stitch are approximately the same and since $2\pi\approx 6$ we get the usual instruction to start with $6$ stitches in a loop and to add $6$ stitches evenly spread on every round. The general formula of adding $2\pi H/W$ stitches allows you to easily check the correct amount of stitches needed for any type of a stitch. We can also calculate how to create a perfect estimate of a sphere or hyperbolic plane using the knowledge of the circumference of circles in spherical and hyperbolic geometry. For more on hyperbolic crochet see \cite{Henderson2001,Taimina2018}.

The plane, sphere, and hyperbolic plane have a constant curvature which means that the circumference of a circle is the same regardless of its centre. However, the crochet models we consider grow in rounds around a fixed centre point and so we only need to be able to calculate the circumference of a circle around this point. An especially interesting candidate for crocheting are the `crochet symmetric' surfaces which can be created with one type of stitch and where the increases or decreases are evenly spaced. This allows simple crochet instructions where only the number of stitches per round needs to be stated. A simple example of such surfaces are the surfaces of revolution. Unfortunately, the catenoid is the only minimal surface that is a surfaces of revolution. The $\B_m$ surfaces are only intrinsically surfaces of revolution but as we will see this is enough to enable simple crochet patterns, at least until intersections. 
  
To find out the required number of stitches per round we need to calculate the circumference of a circle centred at the origin with an intrinsic radius $R=\ell\cdot H$. By intrinsic radius we mean a geodesic, a shortest path between two points, on the surface. 
This process can be split in two steps; 
\begin{itemize}
\item[1.] Calculate the circumference of a circle as a function of the variable $r$ used in \eqref{eq:GeneralEW}.
\item[2.] Determine $r$ as a function of the intrinsic radius $R$. 
\end{itemize}
In the first step we need to determine the arc length of the red curve in Figure \ref{Fig:Surfaces} as a function of $r$ which is given by
\begin{align*}
C(r) = \int_{0}^{2\pi} \sqrt{x_\theta(r,\theta)^2 + y_\theta(r,\theta)^2 + z_\theta(r,\theta)^2} \ d \theta,
\end{align*}
where $x_\theta(r,\theta)=\frac{\partial x(r,\theta)}{\partial \theta}$ is the derivative of $x(r,\theta)$ with respect to $\theta$. From \eqref{eq:GeneralEW} we get
\begin{align*}
x_\theta(r,\theta) &= -r^{m-1}\sin\big((m-1)\theta\big)+r^{m+1}\sin\big((m+1)\theta\big)\\
y_\theta(r,\theta) &= r^{m-1}\cos\big((m-1)\theta\big)+r^{m+1}\cos\big((m+1)\theta\big)\\
z_\theta(r,\theta) &= -2r^m\sin(m\theta).
\end{align*}
Denote $h=h(r,\theta)=(x(r,\theta),y(r,\theta),z(r,\theta))$. Then 
\begin{align*}
\|h_\theta\|_2^2
 &= r^{2(m-1)} + r^{2(m+1)}
 +2r^{2m}\Big(\cos((m-1)\theta)\cos\big((m+1)\theta\big)
 -\sin((m-1)\theta)\sin\big((m+1)\theta\big)\Big)\\
 & \quad +4r^{2m}\sin^2(m\theta)\\
&= r^{2(m-1)} + r^{2(m+1)}+2r^{2m}\cos(2m\theta)+4r^{2m}\sin^2(m\theta)\\
&= (r^{m-1}+r^{m+1})^2,
 \end{align*} 
where we used the trigonometric identities $\cos^2(\theta)+\sin^2(\theta)=1$, $\cos(\theta)\cos(\phi)-\sin(\theta)\sin(\phi)=\cos(\theta+\phi)$, and $\cos(2\theta)=1-2\sin^2(\theta)$. 
We can conclude that the length of the red circle in Figure \ref{Fig:Surfaces} is 
 \begin{align*}
 C(r)
= \int_0^{2\pi}\|h_\theta(r,\theta)\|_2\ d\theta
= 2\pi(r^{m-1}+r^{m+1}).
 \end{align*}
Notice that $\|h_\theta\|_2$ does not depend on the angle $\theta$ which means that the arc length between $\theta_1$ and $\theta_2$ only depends on $|\theta_1-\theta_2|$. This allows us to spread the increases or decreases evenly on every round. 

Since we crochet along the intrinsic radius $R$ instead of $r$ we also need to solve the arc length of the blue curve in Figure \ref{Fig:Surfaces}.
This is given by 
\begin{align}\label{ArcLengthR}
R(\r,\theta) = \int_{0}^{\r} \sqrt{x_r(r,\theta)^2 + y_r(r,\theta)^2 + z_r(r,\theta)^2} \ dr,
\end{align}
where $x_r(r,\theta)=\frac{\partial x(r,\theta)}{\partial r}$ denotes the derivative of $x(r,\theta)$ with respect to $r$. For the Bour's minimal surfaces $\B_m$ we have
\begin{align*}
x_r(r,\theta) &= r^{m-2}\cos((m-1)\theta)-r^{m}\cos\big((m+1)\theta\big)\\
y_r(r,\theta) &= -r^{m-2}\sin((m-1)\theta)-r^{m}\sin\big((m+1)\theta\big)\\
z_r(r,\theta) &= 2r^{m-1}\cos(m\theta).
\end{align*}
Using the same trigonometric identities as before we see that the arc length does not depend on $\theta$; 
\begin{align*}
R(\r,\theta) = R(\r)
&= \int_0^{\r}
r^{m-2}+r^{m}dr\\
&= \frac{\r^{m-1}}{m-1}+\frac{\r^{m+1}}{m+1}.
 \end{align*}
This means that we can create the model using only one type of stitch with height $H$. To find out the number of required stitches on a round $\ell$ we first need to solve $r$ as a function of $R=\ell\cdot H$. Generally this has to be done numerically. We then calculate the circumference of a circle with the intrinsic radius $\ell\cdot H$ and divide this by the width of a stitch $W$, that is, the number of required stitches on the round $\ell$ is $2\pi(r(\ell)^{m-1}+r(\ell)^{m+1})/W$ rounded to the closest integer. 

Above we considered the circumference of circles centred at zero. However, we can infer something more general for the Bour's minimal surfaces. The first fundamental form of a surface is the family of inner products on the tangent spaces induced by the Euclidean inner product on $\R^3$. If we consider the $\B_m$ surface $h(r,\theta)=(x(r,\theta),y(r,\theta),z(r,\theta))$ the coefficients of the first fundamental form are 
\begin{align*}
E_m & =\langle h_r, h_r\rangle = (r^{m-2}+r^m)^2\\
F_m & =\langle h_r,h_\theta\rangle = 0\\
G_m & =\langle h_\theta,h_\theta\rangle = (r^{m-1}+r^{m+1})^2. 
\end{align*}
The first fundamental form completely describes the metric properties of a surface, such as the curvature, length and area. Notice that the first fundamental form of $\B_m$ does not depend on $\theta$. This means that the length of any line segment does not depend on the angle of its position. 
Also, since $F_m=0$ the Gaussian curvature is simply given by 
\begin{align*}
K_m = -\frac{1}{2\sqrt{E_mG_m}}\frac{\partial}{\partial r}\Bigg(\frac{\frac{\partial}{\partial r}G_m}{\sqrt{E_mG_m}}\Bigg)
= -\frac{4r^{2(2-m)}}{(1+r^2)^4}
\end{align*}
and we see that the curvature at a given point only depends on how far the point is from the origin.

\subsection{Enneper surfaces}
\label{Sec:Enneper}

\begin{figure}[H]
\centering
\begin{minipage}[b]{1\textwidth} 
	\includegraphics[width=\textwidth]{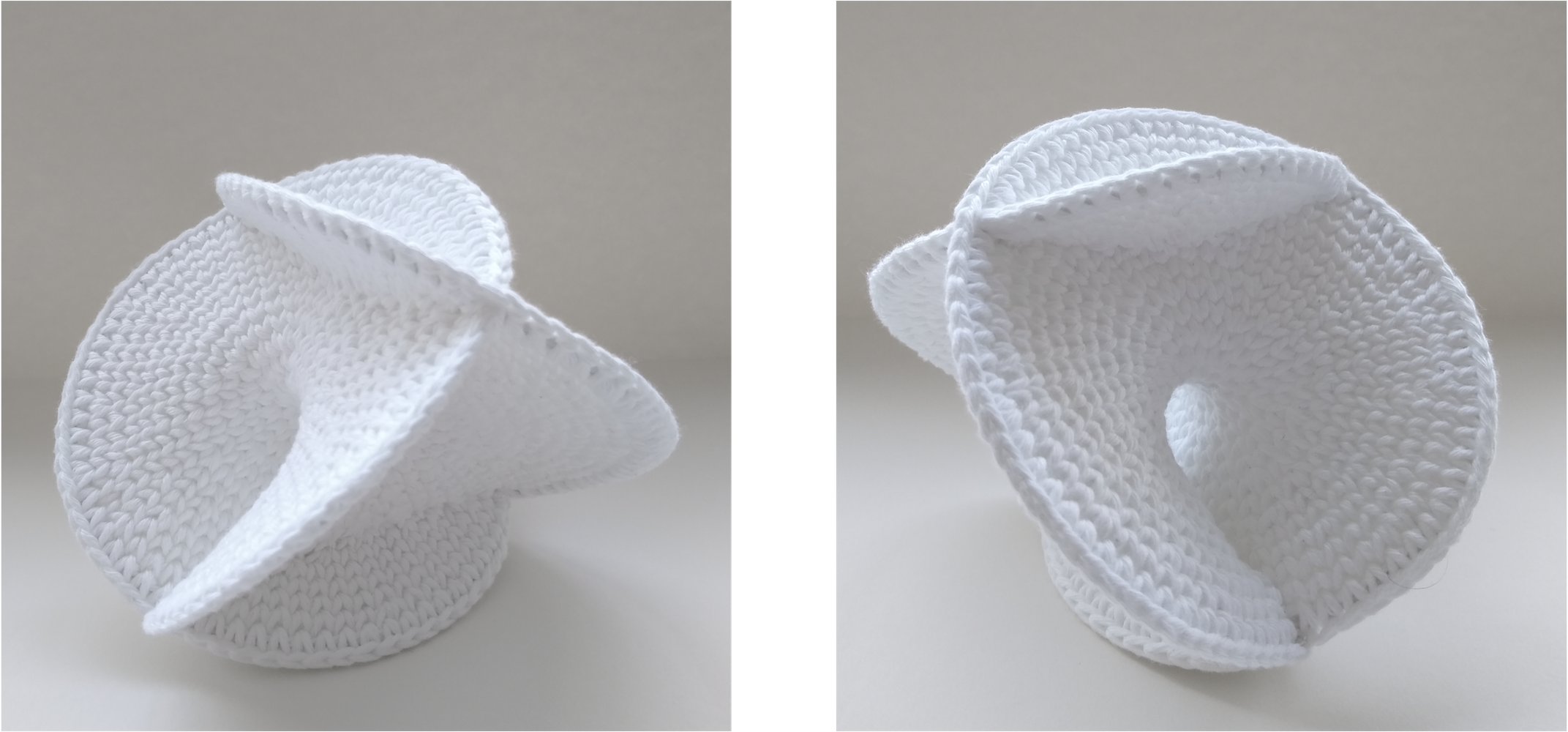}
\end{minipage}
\caption{Two views of the self-intersecting Enneper's surface.}\label{Fig:Enneper1}
\end{figure}

If we choose $m=(k+1)/k$, with $k\in\mathbb{N}_+$, in \eqref{eq:GeneralEW} we attain the family of Enneper's surfaces with order $k+1$ symmetry. See Figure \ref{Fig:Enneper1} for a crocheted Enneper's surface with order 2 symmetry. 
The Enneper's surfaces are often expressed in a simpler form that can be attained with the Weierstrass-Enneper data $f(\zeta)=1$ and $g(\zeta)=\zeta^{n-1}$, $\zeta = \tau e^{i\varphi}$, leading to
\begin{align}\label{eq:Enneper}
\begin{split}
x(\tau,\varphi) &= \tau\cos(\varphi) - \frac{\tau^{2n-1}}{2n-1}\cos\big((2n-1)\varphi\big)\\
y(\tau,\varphi) &= \tau\sin(\varphi) + \frac{\tau^{2n-1}}{2n-1}\sin\big((2n-1)\varphi\big)\\
z(\tau,\varphi) &= \frac{2\tau^n}{n}\cos(n\varphi)
\end{split}
\end{align}
where $\tau\geq0$, $\varphi\in[0,2\pi]$ and $n\in\mathbb{N}$, $n\geq2$, is the order of symmetry. Notice that with $n=1$ the representation \eqref{eq:Enneper} produces a flat disc which is a trivial minimal surface. 
In the following analysis we consider the original parametrisation \eqref{eq:GeneralEW} with $r\geq 0$ and $\theta\in[0,2k\pi]$.

\begin{wrapfigure}{r}{5.5cm}
\includegraphics[width=5.5cm]{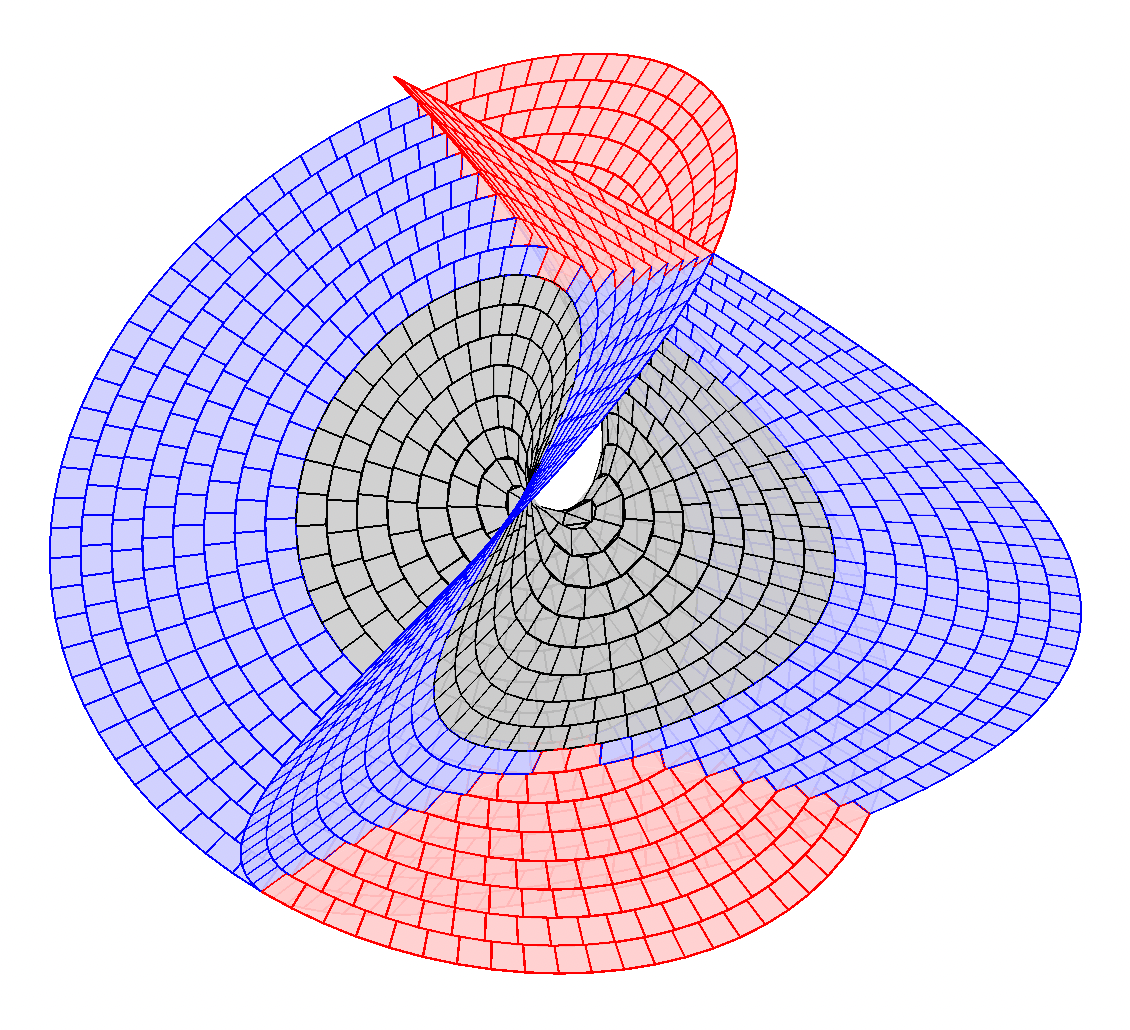}
\caption{Stitches before the intersection in black, stitches in the 'inner' and 'outer' sections after the intersection in red and blue respectively.}
\vspace{-2.5cm}
\label{Fig:EnneperIntersection}
\end{wrapfigure}

For a small $r$ the Enneper surface with order 2 symmetry looks like a saddle but when $r$ grows the surface starts to intersect itself. To crochet a model with the intersections we will have to take a closer look on where the intersection happens at every round. 
The surface first intersects itself at the point $x=y=0$ and $\theta=0$ (or $\theta=\pi$). From this we see that the first intersection occurs when $r=\sqrt{3}$. When $r$ grows we see that the surface has four smaller sections that start to grow, plotted in red in Figure \ref{Fig:EnneperIntersection}, and four larger sections that are getting smaller in comparison, shown in blue. The 
\
crossing angle $\theta_\text{cr}$ satisfies $x(r,\theta)=x(r,\pi-\theta)$, 
that is, 
\begin{alignat*}{2}
        &\;& r\cos(\theta)-\frac{r^3}{3}\cos(3\theta) & = r\cos(\pi-\theta)-\frac{r^3}{3}\cos(3(\pi-\theta))\\
        \Leftrightarrow && 
        \cos(\theta) & = \frac{r^2}{3}(4\cos^3(\theta)-3\cos(\theta))\\
        \Leftrightarrow && 
        \theta & = \cos^{-1}\Bigg(\frac{\sqrt{3}}{2}\sqrt{1+\frac{1}{r^2}}\Bigg),
\end{alignat*}
where we used $\cos(3\theta)=4\cos^3(\theta)-3\cos(\theta)$. Using this, and the fact that the arc length only depends on $|\theta_1-\theta_2|$, we get that the arc length for an `inner' section is $2\theta_\text{cr}(r^{m-1}+r^{m+1})$ and the arc length of an `outer' section is $(\pi/2-2\theta_\text{cr})(r^{m-1}+r^{m+1})$. Notice that the intersection does not happen along a straight line (compare to the $\B_3$ surface in Section \ref{Sec:Bour3}) and some stitches move from the `outer' section into the `inner' section.  This is especially visible right after the intersection occurs. 
The boundary of the `inner' section will become longer than the boundary of the `outer' section when $r>\sqrt{3(1+\sqrt{2})}\approx 2.69$.
We can of course scale \eqref{eq:GeneralEW} to get different sized models requiring a different amount of stitches. 
See Section \ref{Sec:CrochetPractice} for instructions on crocheting your own Enneper surface up to the intersection and Appendix \ref{Appendix} for instructions for a model with intersections.

\subsection{Richmond's surfaces}
The choice of $m = k/(k+1)$, with $k\in\mathbb{N}_+$, in \eqref{eq:GeneralEW} leads to the family of Richmond's surfaces with $k$ order of symmetry. 
The family of Richmond's surfaces can also be expressed by the Weierstrass-Enneper data
$f(\zeta)=\zeta^{-2}$ and $g(\zeta)=\zeta^{k+1}$, $\zeta = \tau e^{i\varphi}$,
which leads to a somewhat cleaner presentation
\begin{align}\label{eq:Richmond}
\begin{split}
x(\tau,\varphi) & = -\frac{\cos(\varphi)}{\tau}-\frac{\tau^{2k+1}}{2k+1}\cos((2k+1)\varphi)\\
y(\tau,\varphi) & = -\frac{\sin(\varphi)}{\tau}-\frac{\tau^{2k+1}}{2k+1}\sin((2k+1)\varphi)\\
z(\tau,\varphi) & = \frac{2\tau^k}{k}\cos(k\varphi),
\end{split}
\end{align}
where $\tau>0$, $\varphi\in[0,2\pi]$, and $k\in\mathbb{N}_+$. We will again use the parametrisation \eqref{eq:GeneralEW}, with $r>0$ and $\theta\in[0,2(k+1)\pi]$, in the analysis below but \eqref{eq:Richmond} is recommendable when plotting the surface. 
Richmond surfaces are sometimes called planar Enneper surfaces because they have one Enneper surface-like self-intersecting end at the middle while the other end stretches out like the plane.
We consider models where $r\in[r_1,r_2]$. Choosing a smaller $r_1$ allows the planar edge to stretch further whereas a larger $r_2$ creates a bigger intersection area in the middle. Notice that the Weierstrass-Enneper data of the Richmond's surfaces with $k=0$ gives a catenoid though it cannot be presented using equations \eqref{eq:GeneralEW} or \eqref{eq:Richmond}.

From a crocheting point of view Richmond's surfaces have two types of self-intersections. We will consider the simplest model with $k=1$. There are two larger `outer' section similar to the Enneper's surface and two smaller `inner' sections, and the surface first intersects when $r=\sqrt{3}$.
The Enneper type intersection angle can be calculated from $y(r,2\pi-\theta)=y(r,2\pi+\theta)$ which leads to
\begin{align*}
\theta_\text{cr}= 2\sin^{-1}\Bigg(\frac{\sqrt{3}}{2}\sqrt{1+\frac{1}{r^2}}\Bigg). 
\end{align*}
Additionally the surface intersects itself on a straight line $x=0$ creating a much more complex crochet pattern. Since the Richmond's surface grows in two direction there is no uniquely best place to start the model. 
Removing crochet stitches usually looks less tidy than adding them so starting from the shortest round and crocheting to the two directions separately is one option.

\subsection{Bour's minimal surface $\B_3$}\label{Sec:Bour3}
The last surface we consider is the Bour's minimal surface $\B_3$ with $r\in[0,1]$ and $\theta\in[0,2\pi]$. 
The Bour's surface $\B_3$ intersects itself on three rays that meet at equal angles at the origin. The rays partition the surface into six equally sized parts. Three of these parts are above the plane spanned by the rays, and three below it. Notice that unlike with the Enneper's and Richmond's surfaces the model consist of six equal parts divided by straight lines and so stitches do not move from one part to another which makes the model easy to crochet after a few set-up rounds.

\begin{figure}[H]
\centering
\begin{minipage}[b]{1\textwidth} 
	\includegraphics[width=\textwidth]{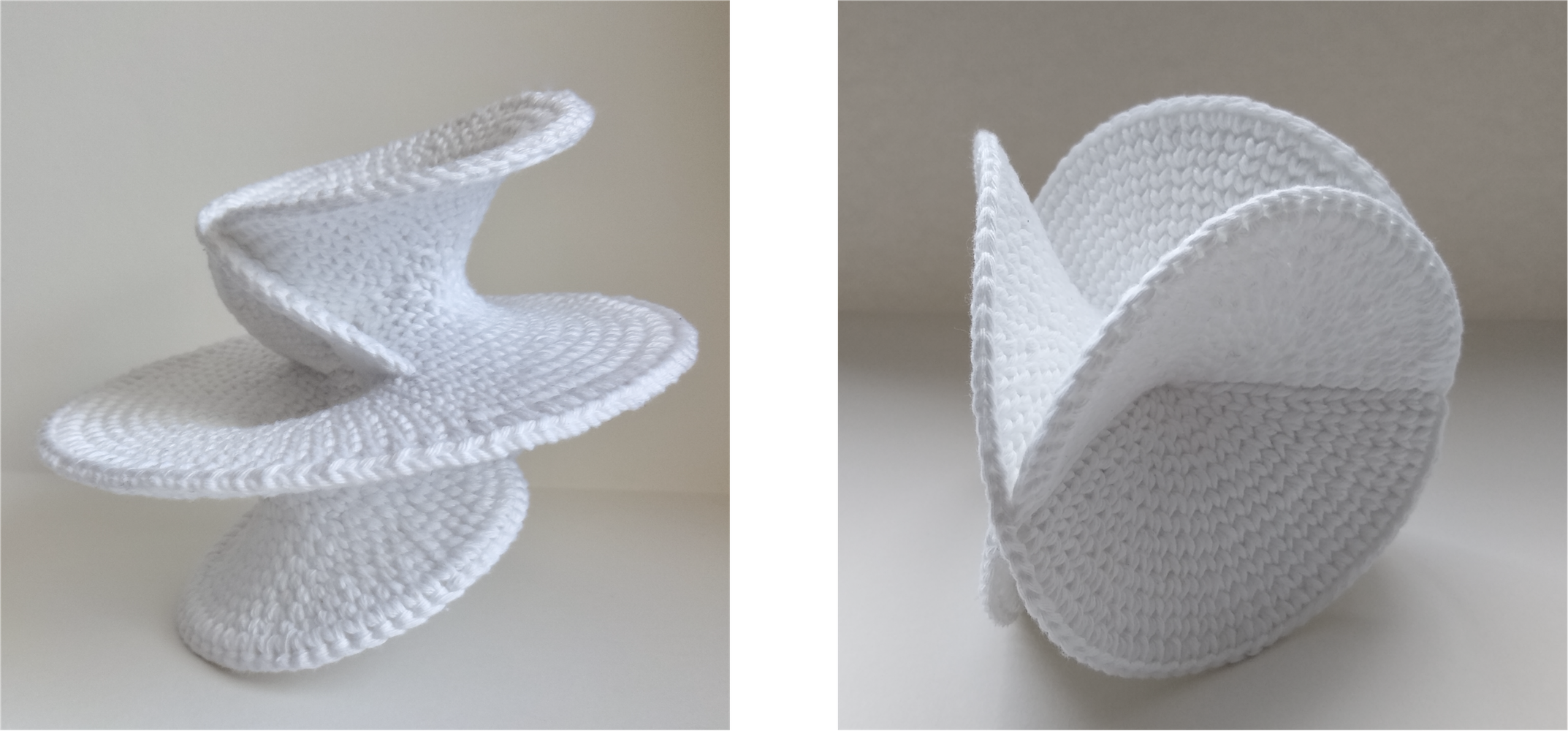}
\end{minipage}
\caption{Crocheted Richmond’s (left) and Bour's (right) minimal surfaces.}\label{Fig:RichmonBour}
\end{figure}

\section{Crochet instructions for the Enneper's surface}\label{Sec:CrochetPractice}

You only need some basic crocheting skills to make your own Enneper's surface. The needed techniques are; magic loop, single stitch, slip stitch, and increase. You can use either the traditional single crochet stitch or the split stitch (aka waistcoat stitch) for a smoother and more rigid surface. If you are a beginner it is recommendable to start with the traditional single crochet stitch since the split stitch is more challenging. You can find many excellent videos online introducing the main crocheting techniques.

\begin{wrapfigure}{r}{5.7cm}
\includegraphics[width=5.7cm]{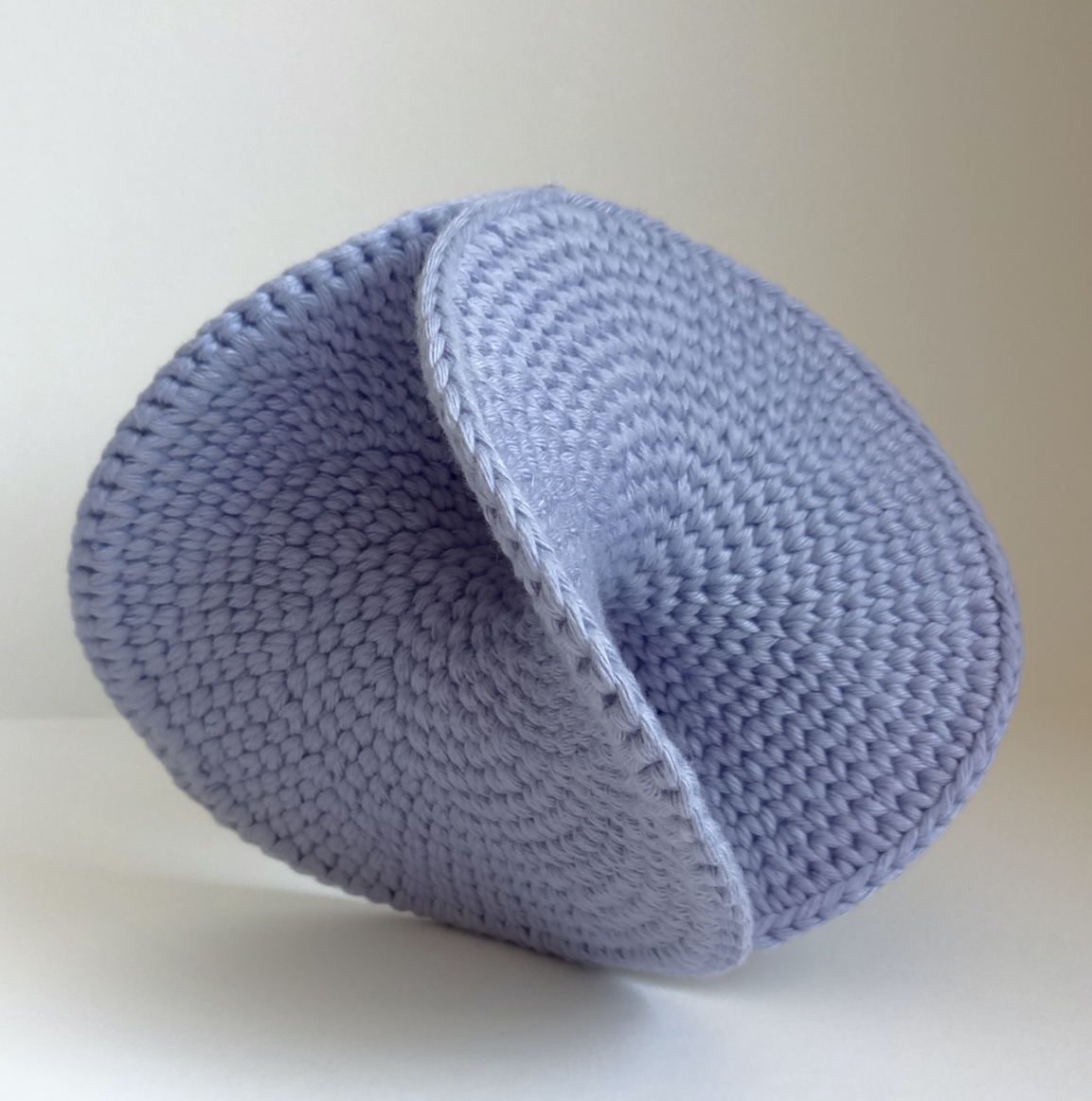}
\caption{Enneper’s surface with 2-fold symmetry crocheted to just before the intersection.}
\label{Fig:EnneperCrocheted}
\end{wrapfigure}

For this project it is best to use a yarn that does not stretch a lot and a crochet hook that is a bit smaller than recommended for the weight to gain a tighter gauge. The models shown in this paper were created using 10-ply (50 grams $\approx$ 75 m) cotton yarn with a 3 mm crochet hook using the split stitch. You should make a test tube to measure the height and width of your stitches. Notice that the model is crocheted in rounds so a tube gives you more accurate   measurements than a traditional swatch. The model is quite sensitive to a wrong stitch size so try to be as consistent as possible. You can find  instructions for three different stitch heights in Table \ref{Tab:Table1}. 
All the models require a bit over $1500$ stitches to complete and the finished model is about $13$ cm wide. Using a thicker (thinner) yarn, while keeping the height width ratio the same, creates a bigger (smaller) model.

Below you can find detailed instructions, in addition to Table \ref{Tab:Table1}, for crocheting the Enneper's surface with gauge W$=0.5$ cm and H$=0.45$ cm up to an intersection as shown in Figure \ref{Fig:EnneperCrocheted}.
The other models are worked in a similar manner. 
To start the model make a magic loop, crochet $6$ stitches into it, and close  the round with a slip stitch. On the next round, $\ell=2$, make two stitches into each stitch from the first round and close the round with a slip stitch. Continue the model according to Table \ref{Tab:Table1} always spacing the indicated number of added stitches evenly on the round. When possible, you should offset your increases so that they are not stacked on top of each other. You can also crochet in continuous rounds, just remember to mark the start of the round. Finishing of the last round will be slightly less tidy when using continuous rounds but the technique does help you to avoid any seams that might appear when closing the rounds with a slip stitch. When you get to the last round you can crochet a thin plastic cable (e.g. transparent trimmer line) into the model. This will help your Enneper's surface to better hold its shape. After you have finished the last round divide the model into four parts and, using a piece of yarn, tie the opposing quarters together making sure that all the four sections have $48$ stitches. Weave off the ends and your model is finished. 

As mentioned above the model is quite sensitive to wrong gauge. If the edge of your model does not curve enough to meet at the last round this is, maybe a bit counterintuitively, because your stitches are too high (or equivalently your stitches are too narrow). Taller stitches mean that on a round $\ell$ the radius of the resulting circle is larger while the circumference stays the same (assuming fixed width) leading to a less curved model. Similarly if it feels that there is no room for the last round this is because your stitches are too short resulting in a too curvy model.

\begin{table}
\begin{subtable}[c]{0.33\textwidth}
\centering
\begin{tabular}{|c|c|c|}
\hline
\multicolumn{3}{|c|}{ H $=0.4$ cm, W $=0.5$ cm}\\
\hline $\ell$ & \hspace*{1.2mm} $\Delta$N \hspace*{1.2mm} & \hspace*{2.5mm} N \hspace*{2.5mm} \\
\hline 1 & & 5 \\
\hline 2 & 6 & 11 \\
\hline 3 & 7 & 18 \\
\hline 4 & 7 & 25 \\
\hline 5 & 9 & 34 \\
\hline 6 & 9 & 43 \\
\hline 7 & 10 & 53 \\
\hline 8 & 10 & 63 \\
\hline 9 & 11 & 74 \\
\hline 10 & 11 & 85 \\
\hline 11 & 11 & 96 \\
\hline 12 & 11 & 107 \\
\hline 13 & 12 & 119 \\
\hline 14 & 12 & 131 \\
\hline 15 & 12 & 143 \\
\hline 16 & 12 & 155 \\
\hline 17 & 12 & 167 \\
\hline 18 & 13 & 180 \\
\hline  &  &  \\
\hline \text {Stitches} & & 1512
\\
\hline
\end{tabular}
\end{subtable}
\begin{subtable}[c]{0.33\textwidth}
\centering
\begin{tabular}{|c|c|c|}
\hline
\multicolumn{3}{|c|}{H $=0.45$ cm, W $=0.5$ cm}\\
\hline $\ell$ & \hspace*{1.2mm} $\Delta$N \hspace*{1.2mm} & \hspace*{2.5mm} N \hspace*{2.5mm}\\
\hline 1 & & 6 \\
\hline 2 & 6 & 12 \\
\hline 3 & 8 & 20 \\
\hline 4 & 9 & 29 \\
\hline 5 & 10 & 39 \\
\hline 6 & 10 & 49 \\
\hline 7 & 12 & 61 \\
\hline 8 & 12 & 73 \\
\hline 9 & 12 & 85 \\
\hline 10 & 12 & 97 \\
\hline 11 & 13 & 110 \\
\hline 12 & 13 & 123 \\
\hline 13 & 13 & 136 \\
\hline 14 & 14 & 150 \\
\hline 15 & 14 & 164 \\
\hline 16 & 14 & 178 \\
\hline 17 & 14 & 192 \\
\hline  &  &  \\
\hline  &  &  \\
\hline \text {Stitches} & & 1525
\\
\hline
\end{tabular}
\end{subtable}
\begin{subtable}[c]{0.33\textwidth}
\centering
\begin{tabular}{|c|c|c|}
\hline
\multicolumn{3}{|c|}{H $=0.5$ cm, W $=0.5$ cm}\\
\hline $\ell$ & \hspace*{1.2mm} $\Delta$N \hspace*{1.2mm} & \hspace*{2.5mm} N \hspace*{2.5mm} \\
\hline 1 & & 6 \\
\hline 2 & 8 & 14 \\
\hline 3 & 9 & 23 \\
\hline 4 & 10 & 33 \\
\hline 5 & 11 & 44 \\
\hline 6 & 12 & 56 \\
\hline 7 & 13 & 69 \\
\hline 8 & 13 & 82 \\
\hline 9 & 14 & 96 \\
\hline 10 & 14 & 110 \\
\hline 11 & 14 & 124 \\
\hline 12 & 15 & 139 \\
\hline 13 & 15 & 154 \\
\hline 14 & 15 & 169 \\
\hline 15 & 15 & 184 \\
\hline 16 & 16 & 200 \\
\hline  &  &  \\
\hline  &  &  \\
\hline  &  &  \\
\hline \text {Stitches} & & 1504
\\
\hline
\end{tabular}
\end{subtable}
\caption{The number added stitches $\Delta$N and the number of stitches N at round $\ell$ for models with gauge W $=0.5$ cm and H $=0.4$ cm, H $=0.45$ cm or H $=0.5$ cm. All the models require a little bit over $1500$ stitches to complete.}
\label{Tab:Table1} 
\end{table}

\newpage
\appendix
\section{Appendix}\label{Appendix}

\begin{wrapfigure}{r}{7cm}
\includegraphics[width=6.5cm]{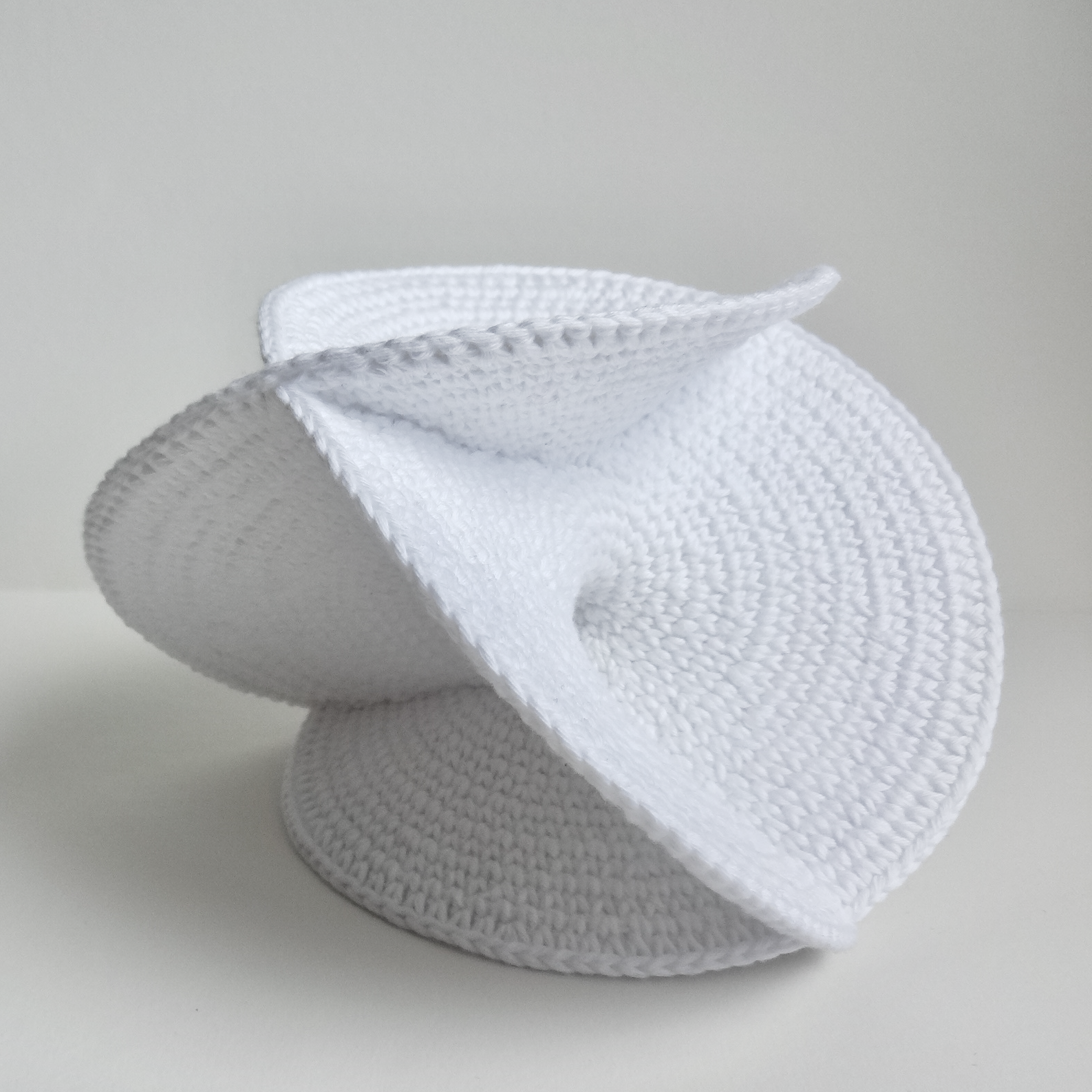}
\caption{Enneper’s surface with a large intersection crocheted with the alternating direction method.}
\label{Fig:EnneperLarge}
\end{wrapfigure}

\begin{figure}[H]
\centering
\begin{minipage}[b]{1\textwidth} 
	\includegraphics[width=\textwidth]{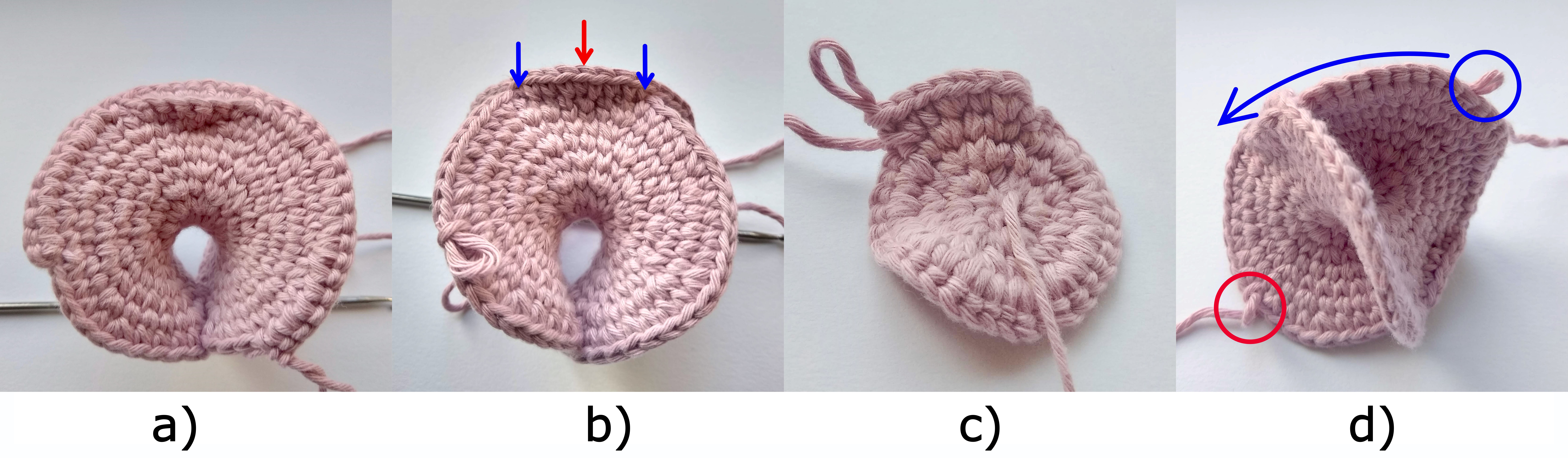}
\end{minipage}
\caption{On the second intersection round crochet two stitches before and after the original intersection stitches to the opposite quarter a). From the third round on you can move one stitch from the outer section into the inner section by skipping the last (first) stitch in the outer section and crocheting an additional stitch into the loop before (after) the first (last) stitch in the inner section, blue arrows in b). Spread the traditional increases evenly from round to round, red arrow in b). To avoid curving of the intersection line caused by tilting of stitches you can alternate the crocheting direction c). If you use the alternating direction method you need to attach a second yarn to the opposite side of the model after the first intersection round d) and complete one round per yarn and then change to the other yarn.}\label{Fig:Appendix2}
\end{figure}

Below you can find instructions on how to crochet a large Enneper's surface with intersections, as shown in Figure \ref{Fig:EnneperLarge}. If you are not familiar with crochet and comfortable with trying some experimental techniques it is recommendable to start with the nonintersecting model described in Section 
An intersecting model requires a few extra considerations. The first one is the intersections themselves. Traditional crochet models do not have intersections and so there are no established methods for doing this. 
The second issue is that crochet stitches are always tilted. This is not a big problem if you are making a model with a small intersection but if you continue the model with evenly spaced increases you will notice that the intersections start to bend. You can try to avoid this by strategically placing increases to counteract the slanting but this technique is somewhat ad hoc and the exact placements depend on how tilted your stitches are. Another way to avoid the unwanted curving of the intersection lines is to alternate the crochet direction, and so the direction of the tilt, from round to round, see Figure \ref{Fig:Appendix2} c). This method will, however, require you to crochet with two balls of yarn after the first intersection round to always cross over one round, see Figure \ref{Fig:Appendix2} d). 

We also need to take into account that the `inner' and `outer' sections do not grow identically. As can be seen in Figure \ref{Fig:EnneperIntersection} some of the stitches move from the outer section into the inner section. For simplicity, the model is split into four equal sized quarters from the last nonintersecting round onwards. We can calculate the arc length of the inner and outer sections as described in Section \ref{Sec:Enneper} and compute the required number of stitches by dividing this by the width of a stitch. After calculating how many increases a quarter has we can calculate how often these increases should appear and from this how many increases an inner and outer section should contain. 
We get the number of stitches moved from an outer section into an inner section as the number of additional stitches per inner section minus the number of increases per inner section, see Table \ref{Tab:Appendix2} for an example.

\begin{table}[htbp!]
\centering
\begin{tabular}{|l|c|c|c|c|c|c|c|c|c|}
 \hline
$\ell$        & 1   & 2   & 3   & 4   & 5   & 6   & 7   & 8   & 9\\
 \hline
 N             & 6   & 14  & 24  & 35  & 46  & 59  & 72  & 86  & 100    \\
 \hline
$\Delta$N   &     & 8   & 10  & 11  & 11  & 13  & 13  & 14   & 14  \\
 \hline
\end{tabular}
\caption{The number of stitches N and the number added stitches $\Delta$N at round $\ell$.}
\label{Tab:Appendix1}   
\end{table}

\begin{figure}[H]
\centering
\begin{minipage}[b]{1\textwidth} 
	\includegraphics[width=\textwidth]{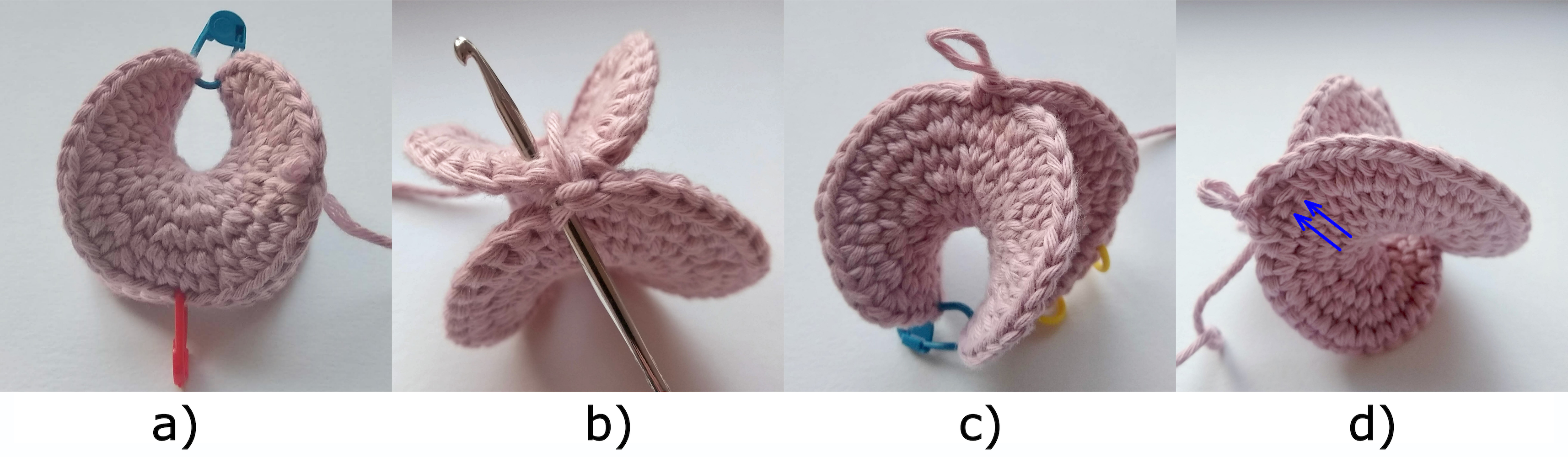}
\end{minipage}
\caption{After you have finished the last round of Table \ref{Tab:Appendix1} divide the model into four equal parts and connect the quarters as shown in a). On the first intersection round crochet the stitches just before and after the first stitch marker to the opposite quarter. To avoid a gap crochet also through the back loop of the original quarter b) and c). This leaves you two stitches (partially) unused d). When you arrive to this intersection again at quarter three jump over and crochet into these unused stitches.}
\label{Fig:Appendix1}
\end{figure}

The model with intersections is calculated for gauge H$=0.45$ cm and W$=0.5$ cm. It is less sensitive for slightly incorrect gauge since the intersection happens quite early. The model is started in the same manner as described in Section \ref{Sec:CrochetPractice} and the number of required stitches per round before intersection can be found in Table \ref{Tab:Appendix1}. After you have finished the last round of the nonintersecting part divide the stitches in four equal parts, each with $25$ stitches, and connect the opposing quarters with stitch markers as shown in Figure \ref{Fig:Appendix1} a). For simplicity, the model is split in four equal parts from now on and the number of stitches below, and in the Table \ref{Tab:Appendix2} refers to the number of stitches per quarter. The inner section, shown in red in Figure \ref{Fig:EnneperIntersection}, will grow both by traditional increases (crocheting twice into a stitch from the previous round) and by moving stitches from the outer section, shown in blue in Figure \ref{Fig:EnneperIntersection}, into the inner section. The different increases, and the increases for the outer section are given in Table \ref{Tab:Appendix2}.

On the first intersection round two stitches will start the inner section. To initiate the intersection crochet until one stitch before the stitch marker and crochet the next two stitches to the opposite quarter as shown in Figure \ref{Fig:Appendix1}. To avoid a gap in the middle crochet also through the back loop of the original quarter, see Figure \ref{Fig:Appendix1} b). On the first intersection round you should make three increases to the outer section.
On the second intersection round the number of stitches in an inner section grows from two to seven. As indicated in Table \ref{Tab:Appendix2}, you should crochet two stitches to the opposite quarter on both sides of the original intersection ($\Delta$N move in $=4$), see also Figure \ref{Fig:Appendix2} a), and make one increase. Make three increases to the outer section. 

On the third round the the number of stitches in an inner section grows from seven to ten by moving two stitches from the outer section into the inner section and by one increase. The simplest way to move one stitch per side is to crochet to the second last stitch on the outer section, jump over the last stitch, and make an extra stitch to the loop before the first stitch in the inner section, see Figure \ref{Fig:Appendix2} b). This loop will most probably be hiding under the stitches from the last round so you will need to push the crossing outwards to reach the loop. When exiting the inner section check very carefully which is the first stitch of the outer section since it will be pulled very tightly towards the inner section and is easy to miss. Make three increases to the outer section.
From now on you will move at most one stitch per side to the inner section and you can use the above method for this. As before you should place the traditional increases evenly and avoid stacking them in the same place from round to round. 

You can finish your model at any round you like. As with the nonintersecting model it is recommendable to crochet a supporting cable into the model on the last round. The full model given in Table \ref{Tab:Appendix2} will be about $20$ cm wide and requires $4394$ stitches. On the last round the length of the inner and outer section is the same requiring $44$ stitches each. 

\begin{table}[htbp!]
\centering
\begin{tabular}{|l|c|c|c|c|c|c|c|c|c|c|c|c|c|c|c|c|c|}
 \hline
$\ell$            & 10  & 11  & 12  & 13  & 14  & 15  & 16  & 17  & 18  & 19  & 20  & 21  & 22  & 23  &   24  &   25  & 26 \\
 \hline
N inner            & 2   & 7   & 10   & 12  & 15  & 17  & 20  & 22  & 25  & 27  & 30  & 32  & 34  & 37  & 39  &  42  &  44\\
 \hline
$\Delta$N move in & 2   & 4   & 2   & 1   & 1   & 1   & 1   & 1   & 1   & 1   & 1  & 1   & 0   & 1   & 0  & 1  & 0\\
 \hline
$\Delta$N inner   & 0   & 1   & 1   & 1   & 2   & 1   & 2   & 1   & 2   & 1   & 2   & 1   & 2   & 2   & 2  & 2  & 2\\
 \hline
 N outer            & 26  & 25  & 26  & 27  & 28  & 29  & 31  & 32  & 33  & 35  & 36  & 37  & 39  & 40  & 42  & 43  & 44\\
 \hline
$\Delta$N outer   & 3   & 3   & 3   & 2   & 2   & 2   & 3   & 2   & 2   & 3   & 2   & 2   & 2   & 2   & 2  & 2  & 1 \\
 \hline
\end{tabular}
\caption{The crochet round $\ell$, the number of stitches in an inner section `N inner', the number of stitches moved from an outer section into an inner section `$\Delta$N move in', the number of traditional increases per inner section `$\Delta$N inner', the number of stitches in an outer section `N outer', and the number of increases per outer section `$\Delta$N outer'.}  
\label{Tab:Appendix2}
\end{table}


\bibliographystyle{aomplain}
\bibliography{ReferencesRecreational}

\providecommand{\bysame}{\leavevmode\hbox to3em{\hrulefill}\thinspace}
\providecommand{\noopsort}[1]{}
\providecommand{\mr}[1]{\href{http://www.ams.org/mathscinet-getitem?mr=#1}{MR~#1}}
\providecommand{\zbl}[1]{\href{http://www.zentralblatt-math.org/zmath/en/search/?q=an:#1}{Zbl~#1}}
\providecommand{\jfm}[1]{\href{http://www.emis.de/cgi-bin/JFM-item?#1}{JFM~#1}}
\providecommand{\arxiv}[1]{\href{http://www.arxiv.org/abs/#1}{arXiv~#1}}
\providecommand{\doi}[1]{\url{https://doi.org/#1}}
\providecommand{\MR}{\relax\ifhmode\unskip\space\fi MR }
\providecommand{\MRhref}[2]{%
  \href{http://www.ams.org/mathscinet-getitem?mr=#1}{#2}
}
\providecommand{\href}[2]{#2}
\begin{thebibliography}{10}

\bibitem{Bour1862}
\bgroup\scshape{}E.~Bour\egroup{}, Theorie de la deformation des surfaces,
  \emph{Journal de l'Ecole Polytechnique} \textbf{22(39)} (1862), 1--148.
  Available at \url{https://gallica.bnf.fr/ark:/12148/bpt6k433694t/f5.item}.

\bibitem{Brown1886}
\bgroup\scshape{}A.~Crum~Brown\egroup{}, On a case of interlacing surfaces,
  \emph{Proceedings of the Royal Society of Edinburgh} \textbf{13} (1886),
  382--386.

\bibitem{Douglas1931}
\bgroup\scshape{}J.~Douglas\egroup{}, Solution of the problem of {P}lateau,
  \emph{Transactions of the American Mathematical Society} \textbf{33} no.~1
  (1931), 263--321. Available at \url{https://doi.org/10.2307/1989472}.

\bibitem{Enneper1868}
\bgroup\scshape{}A.~Enneper\egroup{}, Analytisch-{G}eometrische
  {U}ntersuchungen,  \emph{Nachrichten von der Königl. Gesellschaft der
  Wissenschaften und der Georg-Augusts-Universität zu Göttingen} (1868),
  421--443. Available at \url{http://eudml.org/doc/179395}.

\bibitem{Henderson2001}
\bgroup\scshape{}D.~W. Henderson\egroup{} and
  \bgroup\scshape{}D.~Taimi\c{n}a\egroup{}, Crocheting the hyperbolic plane,
  \emph{Math. Intelligencer} \textbf{23} no.~2 (2001), 17--28. Available at
  \url{https://doi.org/10.1007/BF03026623}.
  
\bibitem{Kekkonen2023}
\bgroup\scshape{}H.N Kekkonen\egroup{}, Calculating crochet models. GitHub repository (2023), Available at
  \url{https://github.com/hnkekkonen/crochet-models}.  

\bibitem{Lagrange1761}
\bgroup\scshape{}J.-L. Lagrange\egroup{}, Essai d'une nouvelle methode pour
  de'terminer les maxima et les minima des formules integrales indefinies,
  \emph{Miscellanea Taurinensia} \textbf{t. II} (1761), 335--362. Available at
  \url{https://gallica.bnf.fr/ark:/12148/bpt6k2155691/f385}.

\bibitem{Meusnier1785}
\bgroup\scshape{}J.~B. Meusnier\egroup{}, M{\'e}moire sur la courbure des
  surfaces,  \emph{Mem des savan etrangers} \textbf{10} (1785 (presented
  1776)), 477--510.

\bibitem{Nitsche1989}
\bgroup\scshape{}J.~C. Nitsche\egroup{}, \emph{Lectures on minimal surfaces:
  vol. 1}, Cambridge university press, 1989.

\bibitem{Osinga2004}
\bgroup\scshape{}H.~M. Osinga\egroup{} and
  \bgroup\scshape{}B.~Krauskopf\egroup{}, Crocheting the {L}orenz manifold,
  \emph{Math. Intelligencer} \textbf{26} no.~4 (2004), 25--37. \mr{2104464}.
  \doi{10.1007/BF02985416}.

\bibitem{Plateau1873}
\bgroup\scshape{}J.~A.~F. Plateau\egroup{}, \emph{Statique exp{\'e}rimentale et
  th{\'e}orique des liquides soumis aux seules forces mol{\'e}culaires},
  \textbf{2}, Gauthier-Villars, 1873.

\bibitem{Rado1930}
\bgroup\scshape{}T.~Rad{\'o}\egroup{}, On {P}lateau's problem,  \emph{Annals of
  Mathematics} (1930), 457--469. Available at
  \url{https://doi.org/10.2307/1968237}.

\bibitem{Reid1971}
\bgroup\scshape{}M.~Reid\egroup{}, The knitting of surfaces,  \emph{Eureka -
  The Journal of the {A}rchimedeans} \textbf{34} (1971), 21--26. Available at
  \url{http://homepages.warwick.ac.uk/~masda/knit_surfaces.pdf}.

\bibitem{Taimina2018}
\bgroup\scshape{}D.~Taimi\c{n}a\egroup{}, \emph{Crocheting adventures with
  hyperbolic planes}, CRC Press, Boca Raton, FL, 2018. Available at
  \url{https://doi.org/10.1201/9780203732731}.

\bibitem{Weierstrass1866}
\bgroup\scshape{}K.~T.~W. Weierstrass\egroup{}, Untersuchungen \"{u}ber die
  fl\"{a}chen, deren mittlere kr\"{u}mmung \"{u}bergleich null ist,  in
  \emph{Mathematische werke von Karl Weierstrass, vol 3}, Berlin Mayer and
  M\"{u}ller, 1894, pp.~39--52. Available at
  \url{https://archive.org/details/mathematischewer03weieuoft/page/n13/mode/2up}.

\bibitem{Whittemore1917}
\bgroup\scshape{}J.~K. Whittemore\egroup{}, Minimal surfaces applicable to
  surfaces of revolution,  \emph{Ann. of Math. (2)} \textbf{19} no.~1 (1917),
  1--20. \mr{1502507}.  \doi{10.2307/1967659}.

\end{thebibliography}

\end{document}